%% file: agt-4-11.tex
\documentclass{gtart_h}

\input agtout

\lognumber{11}
\volumenumber{4}
\volumeyear{2004}
\papernumber{11}
\published{31 March 2004}
\pagenumbers{177}{198}
\received{30 December 2002}
\revised{22 March 2004}
\accepted{25 March 2004}

\usepackage{amssymb}

\input xy
\xyoption{all}

\def\Aut{{\it Aut \,}}
\def\Fix{{\it Fix \,}}

\newtheorem{thm}{Theorem}[section]
\newtheorem{prop}[thm]{Proposition}
\newtheorem{lem}[thm]{Lemma}
\newtheorem{cor}[thm]{Corollary}
\theoremstyle{definition}

\newtheorem{dfn}[thm]{Definition}

\let\demo\proof

\begin{document}

\title{Intersections of automorphism fixed subgroups\\in the free
group of rank three}
\shorttitle{Intersections of automorphism fixed subgroups}
\author{A. Martino} 
\address{Centre de Recerca Matemetica, Apartat 
50\\E-08193 Bellaterra, Spain}
\email{AMartino@crm.es}
\begin{abstract}
We show that in the free group of rank 3, given an arbitrary number of
automorphisms, the intersection of their fixed subgroups is equal to
the fixed subgroup of some other single automorphism.
\end{abstract}

\primaryclass{20E05}\secondaryclass{20F28}
\keywords{Free group, Automorphism, Fixed subgroup}

\maketitle

\section{Introduction}

Let $F_n$ be a free group of rank $n$. 

\begin{dfn}
The \emph{fixed subgroup} of an automorphism $\phi$ of $F_n$,
denoted $\Fix{\phi}$, is the subgroup of elements in $F_n$ fixed
by $\phi$:
 $$
\Fix{\phi}=\{ x \in F_n : \phi x =x \}.
 $$
Following the terminology introduced in~\cite{MaVe}, a subgroup
$H$ of $F_n$ is called \emph{1-auto-fixed} when there exists and
automorphism $\phi$ of $F_n$ such that $H=\Fix{\phi}$. An
\emph{auto-fixed} subgroup of $F_n$ is an arbitrary intersection
of $1$-auto-fixed subgroups. If $S \subseteq \Aut(F_n)$ then
$$
\Fix S=\{ x \in F_n: \phi x =x, \ \forall \phi \in S\}=
\bigcap_{\phi \in S} \Fix{\phi}
$$
is an auto-fixed subgroup. Moreover if $S$ is a subgroup of
$\Aut(F_n)$ generated by $S_0$ then $\Fix S= \Fix S_0$.
\end{dfn}

The celebrated result of Bestvina-Handel~\cite{BH} showed that a
$1$-auto-fixed subgroup of $F_n$ has rank at most $n$. This was
extended by Dicks-Ventura~\cite{DV}. The following Theorem is a special case of the main result in \cite{DV}.
\begin{thm}
\label{inertia}
Let $H$ be a $1$-auto-fixed subgroup of $F_n$ and $K$ a finitely
generated subgroup of $F_n$. Then $rank(H \cap K) \leq rank(K)$. In particular, putting $K=F_n$ we get that $rank(H) \leq n$.  \end{thm}

In~\cite{MaVe} it was conjectured that the families of auto-fixed
and 1-auto-fixed subgroups of $F_n$ coincide. The authors proved
some partial results in this direction but, in general, the
conjecture was only known to be true for $n\leq 2$.

In this paper we show, in Corollary~\ref{rank3}, that any auto-fixed subgroup of $F_3$ is a 
$1$-auto-fixed subgroup of $F_3$. The work in this paper also shows, in Corollary~\ref{upg}, that auto-fixed subgroups are $1$-auto-fixed for an important class of automorphisms called
Unipotent Polynomially Growing (UPG) automorphisms, introduced in \cite{Tits2}.

\medskip

The author gratefully acknowledges the postdoctoral grant SB2001-0128 funded by the Spanish government, and
thanks the CRM for its hospitality during the period of this research.

\section{Preliminaries}

It was shown in \cite{BH} that any outer automorphism
of a finitely generated free group has either polynomial or
exponential growth. We review those notions here.

The growth rate of an automorphism, $\phi$, of a free group $F$ is
the function of $k$ given by the quantity
$$
\sup_{w \in F} \frac{|\phi^k w |}{|w|}
$$
where $|.|$ denotes word length, with respect to a given basis. An
automorphism has polynomial growth if its growth function is
bounded above by a polynomial function of $k$ and exponential
growth if it is bounded below by an exponential function of $k$.
Note that the growth function is always bounded above by an
exponential function of $k$. Clearly, the division into polynomial
and exponential growth automorphisms is independent from the
generating set in question.

One could also replace word length with {\em cyclic} word length -
the word length of a shortest conjugate - and the notions of
polynomial growth and exponential growth are preserved. Thus the
notion of growth rate applies to outer automorphisms, via cyclic
lengths. An outer automorphism, $\Phi$ has polynomial growth if
and only if any automorphism $\phi \in \Phi$ has polynomial
growth. (Having said that, the degrees of the polynomials in
question may differ by at most one.) Similarly $\Phi$ has
exponential growth if and only if $\phi \in \Phi$ has exponential
growth. 

Of particular interest here is the following class, defined in
\cite{Tits2}, Definitions~3.1 and ~3.10.

\begin{dfn}
An outer automorphism, $\Phi$, of the free group of rank $n$,
$F_n$ is said to be Unipotent Polynomially Growing (UPG) if it is
of polynomial growth and its induced action on
$F_n/{F_n'}=\mathbb{Z}^n$, also denoted by $\Phi$,
satisfies one of the following two equivalent conditions.

(i)\qua $\mathbb{Z}^n$ has a basis with respect to which $\Phi$ is
upper triangular with $1$'s on the diagonal, 

(ii)\qua $(Id-\Phi)^k=0$ for some $k>0$.
\end{dfn}

Since we are interested in actual automorphisms and subgroups, we extend the
above notion to these.

\begin{dfn} Let $\Phi \in Out(F_n)$ and let $\phi \in
\Phi$. Then $\phi$ is called UPG if and only if $\Phi$ is UPG. A subgroup of $Out(F_n)$ or $Aut(F_n)$ is called UPG if every element of the subgroup is UPG.
\end{dfn}

A marked graph is a graph, $G$, in the sense of Serre, with a
homotopy equivalence $\tau$ from $R_n$, the rose with $n$ petals
($n$ edges and a single vertex, $*$) to $G$. We identify the free
group of rank $n$, $F_n$, with the fundamental group of $R_n$. The
marking thus gives a specified way to identify the fundamental
group of $G$ with $F_n$.

Suppose that $\tau$ sends $*$ to the vertex $v$ of $G$ and that
$f:G \to G$ is a homotopy equivalence which sends vertices to vertices and edges to edge paths, 
such that $f(v)=u$. If $p$
is any path from $u$ to $v$, then the isomorphism $\gamma_p:
\pi_1(G,u) \to \pi_1(G,v)$ is defined by
$\gamma_p([\alpha])=[\bar{p} \alpha p]$. The map $\gamma_p f_*$ is
then an automorphism of $\pi_1(G,v)$ and the automorphism,
$\phi_{\tau, f,p}$ is then defined so as to make the following
diagram commute.
$$ \xymatrix{ F_n \ar[r]^{\tau_* \ } \ar[dd]_{\phi_{\tau,f,p}} &
\pi_1(G,v) \ar[d]^{f_*} \\ & \pi_1(G,u) \ar[d]^{\gamma_p} \\ F_n
\ar[r]^{\tau_* \ } & \pi_1(G,v) \\
}
$$
As the path $p$ varies amongst all paths between $u$ and $v$, the
collection of automorphisms $\phi_{\tau,f,p}$ form an outer
automorphism $\Phi_{\tau,f}$. 


As in \cite{Tits2}, we shall mostly be interested in the case
where $G$ is filtered, which is to say that there are subgraphs $\emptyset = G_0
\subseteq G_1 \subseteq \ldots \subseteq G_k=G$, where the $G_i$
(not necessarily connected) of $G$. Moreover, each
$G_{i}$ is obtained from $G_{i-1}$ be the addition of a single
edge $E_i$. A map, $f:G \to G$, is {\em upper
triangular with respect to the filtration}, or simply upper
triangular, if for each $i$, $f(E_i)=E_i u_i$, where $u_i$ is a
loop in $G_{i-1}$. Note that this means each edge, and hence the entire graph $G$, has a preferred orientation. It is easy to check that any upper triangular map is a homotopy equivalence on $G$. 

We shall sometimes say that a filtered graph $G$ consists of the edges $E_1, \ldots, E_m$. By this we mean that we have taken one oriented edge from each edge pair of $G$, that the filtration of subgraphs is given by setting $G_i$ to be the subgraph generated by $E_1, \ldots, E_i$ and that the given edges are the preferred orientation for $G$. 

We note that this differs from \cite{Tits2} in that we have
trivial prefixes, which is to say that the image of $E_i$ starts
with $E_i$ rather than another loop in $G_{i-1}$. The more general
situation may be reduced to ours by a sequence of subdivisions. 

For a filtered graph, $G$, one can define $Q$ to be the set of upper triangular homotopy equivalences, 
as above. This turns out to be a group, under composition, which lifts via the marking to a UPG subgroup of $Out(F_n)$. Conversely, a subgroup of $Out(F_n)$ is said to be filtered if it is the lift of a subgroup of such a $Q$ for some $G$.  We rely crucially on the following result from \cite{Tits2}. 

\begin{thm}
\label{filter}
Every finitely generated UPG subgroup of $Out(F_n)$ is filtered.
\end{thm}

We shall also need some further results.

\begin{prop}[\cite{Tits2}, Proposition~4.16]
\label{res}
Let $\phi$ be a UPG automorphism of $F$ and $H\leq F$ a {\em primitive} subgroup of $F$. (That is, $x^k \in H$ implies $x \in H$, for all group elements $x$ and all positive integers $k$.) If $\phi^m(H)=H$ then, in fact, $\phi(H)=H$ and moreover, $\phi$ restricts to a UPG automorphism of $H$. 
\end{prop}

The following shows that periodic points are fixed by UPG automorphisms. 

\begin{prop}
\label{noperiod}
Let $\phi \in \Aut(F_n)$ be a UPG automorphism and suppose that
$\phi^k(x)=x$ for some $x \in F_n$. Then $\phi(x)=x$.
\end{prop}
\demo Without loss we may assume that $x$ is not a proper power
and thus $H=\langle   x\rangle $ is a primitive subgroup. 
By Proposition~\ref{res}, $\phi(H)=H$
and $\phi$ restricts to a UPG outer automorphism of $H$. This implies that $\phi(x)=x$. \endproof

\section{Upper triangular maps}

Throughout this section, $G$ will be a filtered graph with upper triangular maps $f$ and $g$. Note that $f$ and $g$ induce UPG outer automorphisms on the fundamental group of $G$. The edges (and a preferred orientation) of $G$ will be $E_1, \ldots , E_m$ so that the subgraph $G_r$ containing edges $E_1, \ldots, E_r$ 
is invariant under both $f$ and $g$. Recall that $\overline{E}$ denotes the edge with reversed orientation to $E$ and for a path $\alpha$, $\overline{\alpha}$ denotes the inverse path.
By definition we have that, $f(E_i)=E_i u_i$, $g(E_i)=E_i v_i$ where $u_i, v_i$ are loops in $G_{i-1}$. Paths in $G$ will always be edge paths, that is sequences of oriented edges, except for the {\em trivial} paths which will consist of a single vertex. Note that all vertices of $G$ are fixed by both $f$ and $g$. 

We write $\alpha=\beta$ if the paths $\alpha$ and $\beta$ are the same sequence.  We will write $\alpha \simeq \beta $ to denote homotopy equivalence of the paths with respect the endpoints, and write $[\alpha]$ for the homotopy class of the path $\alpha$. 

An edge path is reduced if the sequence of edges contain no adjacent inverse edges and cyclically reduced if the additionally, the fist and last edges are not inverse edges. In each homotopy class of paths, $[\alpha]$ there is a unique reduced path which we denote $\alpha_{\#}$. Write $|\alpha|$ for the number of edges in the path $\alpha_{\#}$. 

Following \cite{Tits1}, section 4.1, the decomposition of a path $\alpha=\alpha_1 \alpha_2 \ldots \alpha_t$ is said to be a splitting for $f$ if for every positive integer $k$, the reduced path $f^k(\alpha)_{\#}$ is obtained by concatenating the reduced paths $f^k(\alpha_i)_{\#}$. Namely, we have that $f^k(\alpha)_{\#} = f^k(\alpha_1)_{\#} \ldots f^k(\alpha_t)_{\#}$. Now, given such a splitting, if we set $\alpha'=\alpha_i \alpha_{i+1}$ then we get another splitting of $\alpha$, $\alpha=\alpha_1 \ldots \alpha_{i-1} \alpha' \alpha_{i+2} \ldots \alpha_t$ which we call a coarsening of the splitting. 

Each path $\alpha$ has a height, denoted $ht(\alpha)$, which is an integer $r$ such that $E_r$ or $\overline{E_r}$ occurs in $\alpha$ but no $E_s$ or $\overline{E_s}$ occurs in $\alpha$ for $s > r$. A {\em basic path of height $r$} is a path of the form $E_r \gamma$, $\gamma \overline{E_r}$ or $E_r \gamma \overline{E_r}$ where $\gamma$ is a path of height less than $r$. The following lemma is proved in \cite{Tits1} for a single map. 

\begin{lem}[\cite{Tits1}, Lemma 4.1.4]
\label{split1}
Let $\alpha$ be a reduced path of height $r$. Then $\alpha$ has a splitting for $f$ into paths which are either basic paths of height $r$ or paths of height less than $r$.
\end{lem}

Note that any path $\alpha$ of height $r$ has a unique decomposition into a minimal number of paths which are either basic of height $r$ or of height less than $r$. Moreover, this unique decomposition is always a coarsening of the decomposition given by Lemma~\ref{split1}. Hence we get the following,

\begin{lem}
\label{split}
Let $\alpha$ be a reduced path of height $r$. Then $\alpha$ has a splitting simultaneously for $f$ and $g$ into paths which are either basic paths of height $r$ or paths of height less than $r$.
\end{lem}

An easy calculation provides the following. 

\begin{lem}
\label{form}
Let $\alpha$ be a basic path of height $r$. If $\alpha$ begins with $E_r$ then so does $f(\alpha)_{\#}$ and if $\alpha$ ends with $\overline{E_r}$ then so does $f(\alpha)_{\#}$.
\end{lem}
\demo The Lemma is clear if $\alpha=E_r \gamma$ or $\gamma \overline{E_r}$. Thus there is only something to show if $\alpha = E_r \gamma \overline{E_r}$. In this case, $\gamma$ is a loop and hence neither $f(\gamma)$ nor $u_r f(\gamma) \overline{u_r}$ can be homotopic to the trivial path and hence $f(\alpha)_{\#}= E_r (u_r f(\gamma) \overline{u_r})_{\#} \overline{E_r}$. \endproof
 
Hence we get the following immediate consequence to Lemmas~\ref{split} and \ref{form}.

\begin{cor}
\label{nontriv}
Let $\alpha$ be a path in $G$ which is not homotopic to a trivial path. Then $f(\alpha)$ and $g(\alpha)$ are not homotopic to trivial paths. 
\end{cor}

\begin{dfn}
Let $\alpha$ be a path of height $r$ which is cyclically reduced. If either the first edge of $\alpha$ is $E_r$ or the last edge is $\overline{E_r}$ then we call $\alpha$ $G$-reduced. 
\end{dfn}

\begin{lem}
\label{stayred}
Let $\alpha$ be a $G$-reduced path. Then $f(\alpha)_{\#}$ is also $G$-reduced. 
\end{lem}
\demo Suppose $\alpha$ begins with $E_r$. We will show that $f(\alpha)_{\#}$ and $g(\alpha)_{\#}$ are cyclically reduced and begin with $E_r$.

First, by coarsening the splitting of Lemma~\ref{split}, there is a splitting of $\alpha=\alpha_1 \alpha_2 \alpha_3$ where $\alpha_1$ is a basic path of height $r$ and $\alpha_3$ is either a basic path of height $r$ or is a path of height less than $r$. Then $\alpha_1$ is either equal to $E_r \gamma$ or $E_r \gamma \overline{E_r}$ for some path $\gamma$ of height less than $r$. Hence, by Lemma~\ref{form}, $f(\alpha_1)_{\#}$ and hence $f(\alpha)_{\#}$ begins with $E_r$. 

So it remains to show that $f(\alpha)_{\#}$ is cyclically reduced, which is to say that it does not end with $\overline{E_r}$. If the path $\alpha_3$ is of height less than $r$ then $f(\alpha_3)$ is not homotopic to a trivial path, by Corollary~\ref{nontriv}, and we are done. Otherwise, $\alpha_3$ is a basic path of height $r$ which cannot end in $\overline{E_r}$ and hence is of the form $E \gamma$ for some path $\gamma$ of height less than $r$. Clearly, $f(\alpha_3)_{\#}$ is not a trivial path and cannot end in $\overline{E_r}$. This proves the Lemma when $\alpha$ begins with $E_r$ and the same proof works for the case when $\alpha$ ends in $\overline{E_r}$ by just repeating the argument for $\overline{\alpha}$.
\endproof

\begin{lem} \label{twored}
Let $\alpha, \beta$ be homotopically non-trivial loops of the same
height based at the same vertex such that $\langle   [\alpha],
[\beta]\rangle $ is a free group of rank $2$ and $\alpha$ is
$G$-reduced. Then there exist positive integers $p,q$ such that
$(\alpha^p \beta \alpha^q)_{\#}$ is $G$-reduced.
\end{lem}
\demo Let $p > |\beta|$. Then, as $\alpha$ is $G$-reduced,
$|\alpha^p|=p|\alpha| \geq p > |\beta|$. Hence, the first letter
of $(\alpha^p \beta)_{\#}$ is the same as the first letter of
$\alpha$. We note that this implies, $$|\alpha^{p+1} \beta|=
|\alpha^p \beta| + |\alpha|.$$
Let $q > |\alpha^{p+1} \beta|> |\alpha^p \beta|$ and consider
$(\alpha^p \beta \alpha^q)_{\#}$. The only way this can fail to be
$G$-reduced is if $(\alpha^p \beta)_{\#} $ is a terminal subpath
of $\overline{\alpha}^q$. Thus $\alpha=\alpha_1 \alpha_2$ (reduced
as written) and
\begin{equation}
\label{eq1} 
\alpha^p \beta \simeq \overline{\alpha_1}
\overline{\alpha^{q_1}}, 
\end{equation}
for some positive integer, $q_1$. Note that we choose $q_1$
maximally so that $\alpha_1 \neq \alpha$.

We repeat this argument for $\alpha^{p+1} \beta \alpha^q$. Either
this is $G$-reduced, and we are done, or $\alpha=\alpha_1'
\alpha_2'$ reduced as written and
\begin{equation} \label{eq2}
\alpha^{p+1} \beta \simeq \overline{\alpha_1'}
\overline{\alpha}^{q_2},
\end{equation}
for some positive integer $q_2$. As before, $\alpha_1' \neq
\alpha$.

The difference in the lengths of the left hand sides of
equations~\ref{eq1} and \ref{eq2} is, as noted above, equal to
$|\alpha|$. Hence the right hand sides must also differ in length
by this amount. Thus,
\begin{eqnarray}\nonumber
|\alpha| & = & |\overline{\alpha_1'} \overline{\alpha}^{q_2}| -
   |\overline{\alpha_1} \overline{\alpha^{q_1}}|  \\\nonumber
&=& (|\alpha_1'| + |\alpha^{q_2}|) - (|\alpha_1| + |\alpha^{q_1}|) \ \mbox{\rm , since $\alpha$ is cyclically reduced} \\
\nonumber & = & (q_2 - q_1)|\alpha| + |\alpha_1'| - |\alpha_1|.
\end{eqnarray}
However, $\alpha_1, \alpha_1'$ are both initial subpaths of
$\alpha$, neither of which is equal to $\alpha$  and so the
quantity $|\alpha_1'| - |\alpha_1|$ must have modulus strictly
less than $|\alpha|$. Hence $\alpha_1 = \alpha_1'$ and
$q_2=q_1+1$. This implies that
$$
\begin{array}{rcl}
\alpha & \simeq & \alpha^{p+1} \beta
\overline{(\alpha^{p} \beta)}, \  \mbox{\rm by inspection}\\
&  \simeq & \overline{\alpha_1} \left.
\overline{\alpha}^{q_1+1}\right. \alpha^{q_1} \alpha_1, \
\mbox{\rm by equations \ref{eq1} and \ref{eq2}}\\ &\simeq&
\overline{\alpha_1} \left. \overline{\alpha} \right. \alpha.\\
\end{array}
$$ This implies that $[\alpha]\neq 1$ is conjugate to its inverse in a free
group. As this cannot happen we get a contradiction and
thus prove the result. \endproof

\begin{prop}
\label{fix} Let $\alpha$ be
a $G$-reduced loop which is fixed up to free homotopy by $f$. Then $f(\alpha) \simeq \alpha$. 
\end{prop}

\demo Now $\alpha$ is a $G$-reduced loop and $f(\alpha)$ is both $G$-reduced by Lemma~\ref{stayred} and freely homotopic to $\alpha$ by hypothesis. However, there are only finitely many $G$-reduced paths freely homotopic to $\alpha$ (they are a subset of the cyclic permutations of $\alpha$). Hence, $f^k(\alpha)_{\#}=\alpha$ for some $k$.  The loop $\alpha$ has a basepoint at some vertex which we denote by $v$. Now $f$ induces a UPG automorphism at $\pi_1(G,v)$ whose $k^{th}$ power fixes $[\alpha]$. By Proposition~\ref{noperiod} we get that $f(\alpha) \simeq \alpha$.  \endproof

The following result is perhaps the key result of the paper. It is an analogue of the Fixed Point Lemma (Corollary~2.2 of \cite{BH}) for a pair of maps. It would be of general interest to see if the following Theorem were true in general, rather than for the UPG case which we restrict ourselves to. 

\begin{thm}
\label{common}
Let $f,g$ be upper triangular maps on $G$. Suppose that $\alpha_1,
\alpha_2$ are loops based at a vertex $v$ in $G$, which generate a
free group of rank $2$ in $\pi_1(G,v)$. Suppose that $\mu, \nu$ are paths
such that $f(\alpha_i) \simeq \mu \alpha_i \overline{\mu}$ and
$g(\alpha_i)\simeq \nu \alpha_i  \overline{\nu}$ for $i=1,2$. Then
there exists a path $\delta$ such that $f(\delta) \simeq \delta
\overline{\mu}$ and $g(\delta) \simeq \delta \overline{\nu}$.
Moreover, $ht(\delta) \leq \mbox{\rm max}\{ht(\alpha_1),
ht(\alpha_2)\}$.
\end{thm}

\demo Consider the set of loops
corresponding to the elements of $\langle    [\alpha_1],
[\alpha_2] \rangle$. To each such loop, $\alpha'$, we can find a $G$-reduced loop, $\alpha$, freely homotopic to it. By Proposition~\ref{fix}, $\alpha$ will be fixed by both $f$ and $g$ up to based homotopy. Amongst all the possible choices, we choose an $\alpha$ which is of minimal height and which is not a proper power. 
Note that there exists a path, which we call $\delta_0$ such that
$\alpha'\simeq \delta_0 \alpha \overline{\delta_0}$. Moreover,
$ht(\delta_0) \leq ht(\alpha') \leq \mbox{\rm max}\{ht(\alpha_1),
ht(\alpha_2)\}$.

Let $\mu_0 \simeq \overline{f(\delta_0)} \mu \delta_0$ and $\nu_0
\simeq \overline{g(\delta_0)} \nu \delta_0$ and note that since
$\alpha$ is fixed up to based homotopy, $\mu_0, \nu_0$ are loops. 
Hence,
$$
\begin{array}{rcl}
\mu_0 \alpha \overline{\mu_0} & \simeq & \overline{f(\delta_0)}
\mu \delta_0 \alpha \overline{\delta_0} \overline{\mu}
f(\delta_0),
\ \mbox{\rm by definition of } \mu_0  \\
& \simeq  & \overline{f(\delta_0)} \mu \alpha'
\overline{\mu} f(\delta_0), \ \mbox{\rm by definition of } \delta_0 \\
& \simeq & \overline{f(\delta_0)} f(\alpha') f(\delta_0), \
\mbox{\rm by hypothesis, since } [\alpha'] \in \langle
[\alpha_1],
[\alpha_2] \rangle \\
& \simeq & \overline{f(\delta_0)} f(\delta_0 \alpha \overline{\delta_0})
f(\delta_0), \ \mbox{\rm by definition of } \delta_0 \\
& \simeq & \alpha, \ \mbox{\rm since } f(\alpha) \simeq \alpha.
\end{array}
$$
As $\alpha$ is not a proper power, we must have that $[\mu_0] \in
\langle    [\alpha] \rangle$. The same calculation for $g$ shows
that $[\nu_0] \in \langle    [\alpha] \rangle$, also.

Choose another loop $\beta'$ representing an element of
$\langle   [\alpha_1], [\alpha_2]\rangle $ so that $[\beta']$ and $[\alpha']$ generate a free group of rank 2.
Now let $\beta =
(\overline{\delta_0} \beta' \delta_0)_{\#}$. 

It is straightforward to verify that $f(\beta) \simeq \mu_0 \beta
\overline{\mu_0}$ and $g(\beta) \simeq \nu_0 \beta
\overline{\nu_0}$. Now, by minimality of $ht(\alpha)$, $ht(\beta) \geq
ht(\alpha)=ht(\mu_0)=ht(\nu_0)$. The proof now separates into two
cases.

\medskip

{\bf Case 1}\qua $r=ht(\beta) > ht(\alpha)$

By Lemma~\ref{split}, $\beta$ has a splitting into paths which are basic paths of height $r$ or paths of height less than $r$. Let $\beta=\beta_1 \ldots \beta_k$ be such a splitting. Note that if $k=1$, then $\beta$ is a basic path of height $r$ and $\beta$ either begins with $E_r$ or ends with $\overline{E_r}$ (or both). But then, by Lemma~\ref{form}, $f(\beta)$ will also either begin with $E_r$ or end with $\overline{E_r}$ and the fact that $f(\beta) \simeq \mu_0 \beta
\overline{\mu_0}$ would mean that $\mu_0$ is homotopic to a trivial path and hence that $f(\overline{\delta_0}) \simeq \overline{\delta_0} \overline{\mu}$. Similarly, $g(\overline{\delta_0}) \simeq \overline{\delta_0} \overline{\nu}$ and we would be finished. Hence, we may assume that $k \geq 2$. 

After coarsening the splitting, we may assume that if $\beta_1$ is a path of height less than $r$, then $\beta_2$ is a basic path of height $r$ starting with $E_r$. Similarly, if $\beta_k$ is a path of height less than $r$, then $\beta_{k-1}$ is a basic path of height $r$ ending in $\overline{E_r}$. Thus, 
$$
f(\beta)_{\#}=f(\beta_1)_{\#} \ldots f(\beta_k)_{\#} \simeq \mu_0 \beta
\overline{\mu_0}.$$
However, the fact that $ht(\mu_0)<r$ implies that $$(\mu_0 \beta
\overline{\mu_0})_{\#}=(\mu_0 \beta_1)_{\#} \beta_2 \ldots \beta_{k-1} (\beta_k
\overline{\mu_0})_{\#}.$$ (Our coarsening of the splitting ensures that $\mu_0$ cannot cause any cancellation with $\beta_2$ or $\beta_{k-1}$). 

Thus we conclude that $f(\beta_k) \simeq \beta_k \overline{\mu_0}$. Hence, 
$$
\begin{array}{rcl}
f(\beta_k \overline{\delta_0}) & \simeq &  \beta_k \overline{\mu_0} f(\delta_0) \\
& \simeq & \beta_k \overline{\delta_0} \overline{\mu}, \ \mbox{\rm  by definition of } \mu_0\\
\end{array}
$$
Similarly, $g(\beta_k \overline{\delta_0}) \simeq \beta_k \overline{\delta_0} \overline{\nu}$ and we would be done in this case.  

\medskip
{\bf Case 2}\qua $ht(\beta)=ht(\alpha)$

By Lemma~\ref{twored} it is possible to find positive integers
$p,q$, such that $\alpha^p \beta \alpha^q$ is $G$-reduced. Then,
by Proposition~\ref{fix}, $\alpha^p \beta \alpha^q$ will be fixed
(up to homotopy rel endpoints) by both $f$ and $g$. However, by
construction
$$\begin{array}{rcl}
f(\alpha^p \beta \alpha^q) &\simeq& \mu_0 \alpha^p \beta \alpha^q
\overline{\mu_0} \\
g(\alpha^p \beta \alpha^q)& \simeq & \nu_0 \alpha^p \beta \alpha^q
\overline{\nu_0}.
\end{array}
$$
Thus,
$$
[\mu_0], [\nu_0] \in \langle   [\alpha]\rangle  \cap \langle
[\alpha^p \beta \alpha^q]\rangle
$$
By definition of $\beta$ this last intersection is trivial, and
hence $\mu_0, \nu_0 \simeq 1$, $f(\delta_0) \simeq \mu \delta_0$,
$g(\delta_0) \simeq \nu \delta_0$. In other words, the theorem is
proved with the requisite $\delta=\overline{\delta_0}$. 
\endproof

The next step is to analyse the fixed paths in $G$. 

\begin{dfn}
A path $\rho$ in $G$ is said to be a common Nielsen Path (NP) for $f$ and $g$ if it is fixed up to homotopy by both $f$ and $g$. $\rho$ is said to be a common Indivisible Nielsen Path (INP) for $f$ and $g$ if it is a common NP and no subpath of $\rho$ is an NP. 
\end{dfn}

The following is immediate from Lemma~\ref{split}.
\begin{lem}
\label{INP}Let $\rho$ be a common INP of height $r$. Then $\rho$ is a basic path of height $r$.
\end{lem}

\begin{prop}
\label{control}  Let $f',g'$ be upper triangular maps on a filtered graph $G'$. Then there is another filtered graph $G=\{E_1, \ldots, E_m\}$, with upper triangular maps $f,g$ and a homotopy equivalence $\tau: G' \to G$ such that the following diagrams commute up to free homotopy. 
$$
\begin{array}{lcr}
\xymatrix{
G' \ar[r]^\tau \ar[d]^{f'} & G \ar[d]^{f} \\
G' \ar[r]^\tau & G\\}
& &
\xymatrix{
G' \ar[r]^\tau \ar[d]^{g'} & G \ar[d]^{g} \\
G' \ar[r]^\tau & G\\
}
\end{array}
$$
Moreover, the maps $f$ and $g$ satisfy the following properties: 

{\rm(i)}\qua If there is a common
Nielsen path of height $r$ then there exist integers $r_f, r_g$
(possibly zero) and a common Nielsen path $\beta_r$ of height at
most $r-1$ such that
$$f(E_r)=E_r {\beta_r}^{r_f}, g(E_r)=E_r {\beta_r}^{r_g}.$$
{\rm(ii)}\qua Up to taking powers, there is a unique common INP of height
$r$. It is $E_r \beta_r \overline{E_r}$ unless $r_f=r_g=0$ in which case it is
$E_r$.
\end{prop}
\demo We will find $G$ from $G'$ by performing a sequence of sliding moves (as in \cite{Tits1}, section~5.4). Now $G'$ has (oriented) edges $E_1', \ldots, E_m'$ so that the subgraph generated by $E_1', \ldots ,E_r'$ is the $r^{th}$ term in the filtration of $G$, and $f(E_r')=E_r' u_r'$, $g(E_r')=E_r' v_r'$. Given any path $\delta$ of height less than $r$ which starts at the terminal vertex of $E_r'$ we can {\em slide $E_r'$ along $\delta$}. Namely, we define a new graph $G$, with the same vertex set as $G'$ and where all the edges $E_i'$ are also edges of $G$ (with the same incidence relations) except for $E_r'$. We have an edge, called $E_r$, of $G$ which starts at the same vertex as $E_r'$ and ends at the same vertex as $\delta$. Informally, we will have that $E_r=E_r' \delta$. The filtration for $G$ will be the same as that for $G'$, just replacing $E_r'$ with $E_r$.  

Define the homotopy equivalence, $\tau$ to be identical on the vertices, and send each $E_i'$ to $E_i'$ for $i \neq r$. Then let $\tau(E_r')= E_r \overline{\delta }$. (Note that $\delta$ is of height less than $r$, so can be considered as a path in both $G$ and $G'$.) 

Now, let $f$ agree with $f'$ on $G'-\{E_r\}$. If the $f'$ image of an edge includes some $E_r'$ we replace each occurrence with $E_r \overline{\delta}$ and then reduce the path. This defines $f$ on all edges except $E_r$ where we let $f(E_r)=E_r (\overline{\delta } u_r' f'(\delta))_{\#}$. Again, note that both $\delta$ and $f'(\delta)$ are paths of height less than $r$. (In fact $f'(\delta)$ is the same path as $f(\delta)$). By construction, $f \tau \simeq \tau f'$. Observe that this is more than an equivalence up to free homotopy, it is also a homotopy equivalence relative to the vertex sets of $G$ and $G'$. The inverse map to $\tau$ is one which sends $E_i'$ to $E_i'$ for $i \neq r$ and sends $E_r$ to $E_r' \delta$. Similarly, we can define $g$ and note that $g(E_r)= E_r (\overline{\delta } v_r' g'(\delta))_{\#}$.

Now the fact that $\tau$ is a homotopy equivalence relative to the vertex sets means that if $f'(E_i')=E_i' u_i'$ for $i \neq r$ where $u_i'$ is a NP then $\tau(u_i')$ is also a NP path for $f$ and $f(E_i')=E_i' \tau(u_i')_{\#}$. Also it is clear that both $\tau$ and its homotopy inverse preserve the height of paths and send basic paths to basic paths of the same type. The strategy is to perform sliding moves so as to make properties (i) and (ii) hold for as many edges as possible. The above comments show that if we slide $E_r'$ along some path, we do not disturb properties (i) and (ii) for other edges. 

Thus it is sufficient to show that we may perform sliding homotopies for the edge $E_r'$ so that properties (i) and (ii) hold for the resulting maps $f$ and $g$ with respect to $E_r$.

Suppose there is a common Nielsen path for $f'$ and $g'$ of height $r$. Hence there must be
a common indivisible Nielsen path of height $r$ which, by Lemma~\ref{INP} and up to orientation, 
must be of the form $E_r' \gamma$ or $E_r' \gamma \overline{E_r'}$
where $\gamma$ is a path of height at most $r-1$. In the former
case we slide $E_r'$ along $\gamma$ and then both $f$ and $g$ fix $E_r$, to satisfy (i) and (ii). 

So we shall assume that there is no common INP of the form $E_r'
\gamma$. Thus we are in the latter case, and there is a common INP of
the form $E_r' \gamma \overline{E_r'}$, where $\gamma$ is not a proper power. 

Now we can find a path, $\delta$, of height less than $r$ such that $(\overline{\delta} \gamma {\delta})_{\#}$ is $G'$ reduced. Let $\beta_r=(\overline{\delta} \gamma {\delta})_{\#}$. Then by, Proposition~\ref{fix}, $\beta_r$ is fixed by $f'$ and $g'$. Moreover, 
$$
\begin{array}{rcl}
\beta_r & \simeq & f'(\beta_r), \mbox{\rm as $\beta_r$ is fixed} \\
& \simeq & f'(\overline{\delta} \gamma {\delta}) \\
& \simeq & f'(\overline{\delta}) f'(\gamma) f'(\delta) \\
& \simeq & f'(\overline{\delta})\overline{u_r'} \gamma u_r'  f'(\delta), \mbox{\rm as $E_r' \gamma \overline{E_r'}$ is fixed} \\
& \simeq & f'(\overline{\delta}) \overline{u_r'} \delta \beta_r \overline{\delta} u_r'  f'(\delta)\\
\end{array}
$$
Hence, $(\overline{\delta} u_r'  f'(\delta))_{\#}=\beta_r^{r_f}$ for some integer, $r_f$. Hence, if we slide $E_r'$ along $\delta$, we get that $f(E_r)=E_r \beta_r^{r_f}$, where $\beta_r$ is a NP of both $f$ and $g$ of height less than $r$. Similarly, $g(E_r)=E_r \beta_r^{r_g}$ for some integer, $r_g$. 

Now if for our new maps, $f$, $g$, the NP $E_r \beta_r \overline{E_r}$ is the only NP which is a basic path of height $r$, then it must be an INP and we would be done. By assumption, there is no common NP of the form $E_r \gamma$ (strictly speaking, we have assumed there is no such in $G'$, but the existence of such a common NP for $f$ and $g$ in $G$ implies that there is also one in $G'$ for $f'$ and $g'$). 

So let us assume that there is a common NP of the form $E_r \alpha \overline{E_r}$ for $f$ and $g$, where $\alpha$ is not a power of $\beta_r$, and arrive at a contradiction. 

Then, $[\alpha], [\beta]$ are loops based at the same vertex and must generate a free group of rank
$2$. Moreover, $f(\alpha) \simeq \overline{\beta_r^{r_f}} \alpha \beta_r^{r_f}$ and $g(\alpha) \simeq \overline{\beta_r^{r_g}} \alpha \beta_r^{r_g}$ and recall that $\beta_r$ is a NP for both $f$ and $g$. 
Hence, by Proposition~\ref{common}, there exists a path $\delta_0$ of height less than $r$ such that $f(\delta_0) \simeq \delta_0 \beta_r^{r_f}$ and $g(\delta_0) \simeq \delta_0 \beta_r^{r_g}$. This means that $E_r \overline{\delta_0}$ is a common NP for $f$ and $g$. This contradiction completes the proof. \endproof

\begin{prop}
\label{allfix} Let $f,g$ be Upper triangular maps satisfying the
conclusion of Proposition~\ref{control}. If there is no common INP
of height $r$ then no common Nielsen path can cross $E_r$.
\end{prop}
\demo We argue by contradiction. Let $\rho$ be a common  Nielsen
path of height $k \geq r$ which crosses $E_r$ and choose $k$ to be
the smallest integer with this property. The point is that $E_r$
is an edge in $\rho$ but, a priori, we do not know if there is an
INP of height $r$ as a subpath of $\rho$.

We decompose $\rho$ into common INP's. If $E_k$ is a common INP
itself then this decomposition is into paths each of which is
equal to $E_k$ or is a common INP of height $\leq k-1$. Thus
either $r=k$ or $E_r$ is an edge in a in a common INP of height at
most $k-1$. The former contradicts the hypotheses and the latter
the minimality of $k$.

If $E_k$ is not a common INP then by Proposition~\ref{control},
every common INP of height $k$ is equal to $E_k \beta_k^m
\overline{E_k}$ for some integer $m$,  where $\beta_k$ is a common
INP of height at most $k-1$.
Thus $\rho$ can be written as a product of common INP's
each of which is equal to (a power of) $E_k \beta_k
\overline{E_k}$ or is an INP of height at most $k-1$. Thus either
$E_r=E_k$, contradicting the hypotheses or $E_r$ is an edge in a
Nielsen path of height at most $k-1$, contradicting the minimality
of $k$. \endproof

\begin{thm}
\label{abelian} Let $\phi, \psi \in Aut(F_n)$ such that $\langle
\phi, \psi\rangle $ is a UPG subgroup and $\Fix \phi \cap
\Fix \psi$ is a free group of rank at least two which is contained
in no proper free factor of $F_n$. Then there exists a filtered graph
$G=\{E_1,...,E_m\}$, upper triangular maps $f,g:G \to G$ and an isomorphism $\tau:F_n \to \pi_1(G,v)$ such that the following diagrams commute:
$$
\begin{array}{lcr}
\xymatrix{
F_n \ar[r]^{\tau \ \ } \ar[d]^\phi & \pi_1(G,v) \ar[d]^{f_*} \\
F_n \ar[r]^{\tau \ \ }  & \pi_1(G,v)\\}
& &
\xymatrix{
F_n \ar[r]^{\tau \ \ } \ar[d]^\psi & \pi_1(G,v) \ar[d]^{g_*} \\
F_n \ar[r]^{\tau \ \ } & \pi_1(G,v).\\
}
\end{array}
$$
Moreover, $f(E_r)=E_r
{\beta_r}^{r_f}$ and $g(E_r)=E_r {\beta_r}^{r_g}$ where $\beta_r$ is a
common NP of height at most $r-1$ and $r_f, r_g $ are integers. The only common INP of height $r$ is 
(a power of) $E_r {\beta_r} \overline{E_r}$, unless $r_f=r_g=0$, in which case it is $E_r$. 
\end{thm}

\demo By Theorem~\ref{filter} we can find a upper triangular maps $f,g:G \to G$,
which represent the {\em outer} automorphisms corresponding to
$\phi$ and $\psi$. Thus we have that the following diagrams commute up to free homotopy. 
$$
\begin{array}{lcr}
\xymatrix{
R_n \ar[r]^\tau \ar[d]^\phi & G \ar[d]^f \\
R_n \ar[r]^\tau & G\\} & & \xymatrix{
R_n \ar[r]^\tau \ar[d]^\psi & G \ar[d]^g \\
R_n \ar[r]^\tau & G\\
}
\end{array}
$$
Here $\tau$ is a homotopy equivalence and $R_n$ is the rose with $n$ edges. By Proposition~\ref{control}, we can assume that $f$ and $g$ satisfy the conclusions of that result. Next, we know that we can find paths $\mu, \nu$ such that the following diagrams commute
$$
\begin{array}{lcr}
\xymatrix{ F_n \ar[r]^{\tau_* \ } \ar[dd]_\phi &
\pi_1(G,v) \ar[d]^{f_*} \\ & \pi_1(G,v) \ar[d]^{\gamma_{\mu}} \\ F_n
\ar[r]^{\tau_* \ } & \pi_1(G,v) \\
}& & \xymatrix{ F_n \ar[r]^{\tau_* \ } \ar[dd]_\psi &
\pi_1(G,v) \ar[d]^{g_*} \\ & \pi_1(G,v) \ar[d]^{\gamma_{\nu}} \\ F_n
\ar[r]^{\tau_* \ } & \pi_1(G,v) \\
}
\end{array}
$$
Here, for a loop $p$, $\gamma_p$ is the inner automorphism induced by $p$, $\gamma_p(\alpha) \simeq \overline{p} \alpha p$. 

As $\Fix \phi \cap \Fix \psi$ has rank at least two, we may apply
Theorem~\ref{common} to find a path $\delta$ such that $f(\delta)
\simeq \delta \mu$ and $g(\delta) \simeq \delta
\nu$. Hence we get the following commuting diagrams. 
$$
\begin{array}{lcr}
\xymatrix{
 \pi_1(G,v) \ar[r]^{\gamma_{\delta}^{-1}} \ar[d]^{f_*} & \pi_1(G,v) \ar[dd]^{f_*} \\ 
\pi_1(G,v) \ar[d]^{\gamma_{\mu}} &  \\
\pi_1(G,v)  \ar[r]^{\gamma_{\delta}^{-1}} & \pi_1(G,v) \\
}  
& 
\xymatrix{
 \pi_1(G,v) \ar[r]^{\gamma_{\delta}^{-1}} \ar[d]^{g_*} & \pi_1(G,v) \ar[dd]^{g_*} \\ 
\pi_1(G,v) \ar[d]^{\gamma_{\nu}} &  \\
\pi_1(G,v)  \ar[r]^{\gamma_{\delta}^{-1}} & \pi_1(G,v) \\
}
\end{array}
$$
Thus we get a commuting diagram as required by the statement of this Theorem and, without loss of generality, we may also assume that there are no valence one vertices in $G$. If, for some $r$, there is no common INP of height $r$ then we may apply Proposition~\ref{allfix} to deduce that the image of $\Fix \phi \cap \Fix \psi$ in $\pi_1(G,v)$ is a subgroup generated by loops, none of which cross $E_r$. This would imply that $\Fix \phi \cap \Fix \psi$ is contained in a proper free
factor, as $G$ has no valence one vertices, and hence would be a contradiction. Thus there must be an INP at every height and we are done since we already know that $f$ and $g$ satisfy the conclusion of Proposition~\ref{control}. \endproof

In order to prove that the fixed subgroup of any UPG subgroup is 1-auto-fixed, we need to invoke the following 
result.

\begin{thm}[\cite{MV1}, Corollary~4.2] 
\label{finite}
For any set $S \subseteq Aut(F_n)$, there exists a finite subset $S_0 \subseteq S$ such that $\Fix(S)=\Fix(S_0)$.
\end{thm}

\begin{cor}
\label{noff} Let $A \leq \Aut{F_n}$ be a UPG
subgroup such that $\Fix{A}$ is contained in no proper free
factor. Then there exists a $\chi \in \Aut(F_n)$ such that
$\Fix{A}=\Fix{\chi}$. If the rank of $\Fix{A}$ is at least $2$, then we may choose $\chi \in A$.
\end{cor}
\demo If $\Fix{A}$ is cyclic, generated by $w$, then $w$ 
cannot be a proper power and we can take
$\chi$ to be the inner automorphism by $w$. 

Otherwise, by Theorem~\ref{finite}, we note that it 
is sufficient to prove this Corollary for the case where
$A$ is generated by two elements, $\phi, \psi$. Apply
Theorem~\ref{abelian} and, using the notation from that Theorem,
consider the Upper triangular map $fg^k$ for some integer $k$.
Clearly, $fg^k(E_r)=E_r \beta_r^{r_f+k.r_g}$. If $k$ is
sufficiently large then $r_f+k.r_g$ is only equal to $0$ if
$r_f=r_g=0$. For this value of $k$, a path is fixed by $fg^k$ if
and only if it is fixed by both $f$ and $g$. Hence $\Fix{\phi
\psi^k}=\Fix( \langle \phi,\psi \rangle)$. \endproof

\begin{cor}
\label{upg}
Let $A \leq \Aut{F_n}$ be a UPG subgroup. Then
there exists a $\chi \in \Aut{F_n}$ such that $\Fix{\chi}=\Fix{A}$
\end{cor}

\demo Let $H$ be a free factor of $F_n$ of smallest rank
containing $\Fix{A}$. Let $\phi \in A$, then $\phi(H)$ is another
free factor of $F_n$ containing $\Fix{A}$. Hence $H \cap \phi(H)$
is a free factor of $F_n$ and hence also a free factor of both $H$
and $\phi(H)$. By minimality of the rank of $H$, $H \cap \phi(H)$
must be equal to $H$ and hence $H$ is a free factor of $\phi(H)$.
As the rank of $H$ is equal to that of $\phi(H)$, $H=\phi(H)$ for all $\phi \in A$.

Thus, we may look at the restriction of $A$ to $Aut(H)$. By \cite{Tits2} Proposition~4.16, 
this is a UPG subgroup. Clearly, $\Fix{A}$ is contained is
no proper free factor of $H$ and hence by Corollary~\ref{noff},
there exists a $\widehat{\chi} \in Aut(H)$ such that
$\Fix{\widehat{\chi}}=\Fix{A}$.

As $H$ is a free factor of $F_n$, we may find a basis $x_1, \ldots
,x_n$ of $F_n$ such that $H=\langle    x_1, \ldots ,x_r \rangle$
for some $r$. Define $\chi \in Aut(F_n)$ by letting the $\chi$
agree with $\widehat{\chi}$ on $H$ and letting
$\chi(x_j)={x_j}^{-1}$ for $j > r$. Clearly,
$\Fix{\chi}=\Fix{\widehat{\chi}}$, and we are done. \endproof

\section{The rank $n-1$ case}

In this section we describe the structure of fixed subgroups of
exponential automorphisms where the fixed subgroup has rank one
less than the ambient free group. In order to do this, we invoke
the main Theorem of \cite{MV}, but first a definition.

\begin{dfn}
Let $\phi$ be an automorphism of $F$ and $H$ a subgroup of $F$. $H$ is called $\phi$-invariant if $\phi(H)=H$ setwise.  
\end{dfn}

\begin{lem}
\label{poly} Let $F_n$ be a free group with basis $x_1, \ldots
,x_n$ and let $H$ be the free factor $\langle    x_1, \ldots, x_{n-1} \rangle$. Suppose that for some
$\phi \in Aut(F_n)$, $H$ is $\phi$-invariant. Then $\phi(x_n)=u
{x_n}^{\pm 1} v$ for some $u,v \in H$. Moreover, if $\phi|_H$ has
polynomial growth, then so does $\phi$.
\end{lem}

\demo As $\phi$ is an automorphism, $\phi(x_1), \ldots, \phi(x_n)$
is a basis for $F_n$. By applying Nielsen moves only on $H$ we may
deduce that $x_1, \ldots, x_{n-1}, \phi(x_n)$ is also a basis for
$F_n$. (Alternatively, one can extend $\phi|_H$ to an
automorphism, $\phi'$ of $F_n$ by letting $\phi'(x_n)=x_n$. The
image of $x_1, \ldots, x_n$ under the automorphism $\phi
{\phi'}^{-1}$ shows that $x_1, \ldots, x_{n-1}, \phi(x_n)$ is a
basis.)

Now we write $\phi(x_n)=uwv$ where this product is reduced as
written, $u,v \in H$ and $w$ is a word whose first and last
letters are ${x_n}^{\pm 1}$. Clearly, the set $x_1, \ldots,
x_{n-1}, w$ is a basis. Moreover it is a Nielsen reduced basis
(see \cite{LS}) and hence $w={x_n}^{\pm 1}$. This proves the first
claim.

To prove the second claim we note that if $\phi|_H$ has polynomial
growth then $|\phi^k(g)| \leq |g| A k^d$ for all $g \in H$ and
some constants $A,d$. Hence,
$$
\begin{array}{rcl}
|\phi^k(x_n)| & \leq & \sum_{i=0}^{k-1}
|\phi^i(u)|+\sum_{i=0}^{k-1} |\phi^i(v)|+1 \\
& \leq & (|u|+|v|+1) \sum_{i=0}^{k-1} A i^d\\
& \leq & (|u|+|v|+1)A k^{d+1},\\
\end{array}
$$
from which it follows easily that $\phi$ also has polynomial
growth (where the degree of the polynomial is at most one higher).
\endproof

The following result gives a description of 1-auto-fixed subgroups which will be particularly useful in our situation.

\begin{thm}{\rm\cite{MV}}\label{description}\qua
Let $F$ be a non-trivial finitely generated free group\break and let
$\phi\in Aut(F_n)$ with $\Fix{\phi}\neq 1$. Then, there exist
integers $r\geq 1$, $s\geq 0$, $\phi$-invariant non-trivial
subgroups $K_1, \ldots ,K_r \leq F_n$, primitive elements $y_1,
\ldots ,y_s \in F_n$, a subgroup $L\leq F_n$, and elements $1\neq
h'_j \in H_j =K_1*\cdots *K_r*\langle    y_1,\ldots ,y_j \rangle$,
$j=0,\ldots ,s-1$, such that
 $$
F_n=K_1*\cdots *K_r*\langle    y_1,\ldots ,y_s \rangle *L
 $$
and $\psi y_j =y_j h'_{j-1}$ for $j=1,\ldots ,s$; moreover,
 $$
\Fix{\phi}=\langle    w_1, \ldots ,w_r, y_1 h_0 y_1^{-1}, \ldots
,y_s h_{s-1}y_s^{-1} \rangle
 $$
for some non-proper powers $1\neq w_i\in K_i$ and $1\neq h_j\in
H_j$ such that $\phi h_j =h_j' h_j h_j^{\prime -1}$, $i=1,\ldots
,r$, $j=0,\ldots ,s-1$. \qed
\end{thm}

\begin{prop}{\rm\cite{CT}}\qua
\label{maxpoly} Let $\phi \in Aut(F_n)$ and suppose that
$rank(\Fix{\phi})=n$. Then $\phi$ is UPG.
\end{prop}
\demo This follows directly from \cite{CT}, although it is not
stated in these terms. It can also be proven using
Theorem~\ref{description}, since the hypothesis guarantees that $L$
is trivial and each $K_i$ is cyclic. Since it is clear that $H_j$ is $\phi$-invariant, 
repeated applications of Lemma~\ref{poly} show
that $\phi$ has polynomial growth. We can then deduce that $\phi$
is UPG by inspection of the basis given in
Theorem~\ref{description}. \endproof

\begin{prop}
\label{basis} Let $\chi \in \Aut{F_n}$ be an automorphism of
exponential growth and suppose that $rank(\Fix{\chi})=n-1$. Then there
exists a basis
$x,y,a_1,...,a_{n-2}$ of $F_n$ such that

{\rm(1)}\qua Each free factor, $A_i=\langle x,y,a_1, \ldots ,a_i \rangle$, is $\chi$-invariant. In particular, 
$\langle x,y \rangle$ is $\chi$-invariant and $\chi|_{\langle x,y \rangle}$ has exponential growth.

{\rm(2)}\qua  $\chi a_i = a_i w_i$ where $w_i \in \langle x,y,a_1,...,a_{i-1} \rangle$.

{\rm(3)}\qua  $\Fix{\chi}=\langle xyx^{-1}y^{-1}, y_1,...,y_{n-2}\rangle$
where each $y_i$ is equal to either $a_i$ or $a_i g_i a_i^{-1}$
for some $g_i \in \langle x,y,a_1,...,a_{i-1} \rangle$.
\end{prop}

\demo We apply Theorem~\ref{description}. First note that
$$F_n=K_1*\cdots *K_r*\langle    y_1,\ldots ,y_s \rangle
*L=H_s*L$$ and that $H_s$ is $\chi$ invariant. Moreover,
$n-1=r(\Fix{\chi})=r+s \leq r(H_s)$.

Hence if $L$ is not trivial, then $n-1=rank(\Fix{\chi}) \leq rank(H_s)
\leq n-1$. Thus $rank(\Fix{\chi})=H_s$ and $rank(L)=1$. Thus, by
Lemma~\ref{maxpoly}, $\chi|_{H_s}$ has polynomial growth and we
deduce that $\chi$ has polynomial growth by Lemma~\ref{poly}. This
is a contradiction as we are assuming that $\chi$ has exponential
growth, by hypothesis. Hence $L=1$.

Since $r+s=n-1$ and $n=rank(H_s)=s+\sum_{i=1}^r rank(K_i)$ it is easy to
see that all but one of the $K_i$ is cyclic and that this exception
must have rank 2. Without loss of generality, $rank(K_1)=2$ and every
other $K_i$ is cyclic. Let $x,y$ be a basis for $K_1$ and extend
this to a basis $x,y,a_1, \ldots,a_{n-2}$ where $a_1,
\ldots,a_{r-1}$ are the generators of $K_2, \ldots, K_r$
respectively and $a_r, \ldots , a_{n-2}$ are equal to $y_1,
\ldots, y_s$ respectively.

By construction, $\chi a_i = a_i w_i$ where $w_i \in \langle
x,y,a_1,...,a_{i-1} \rangle$. Moreover,
$$\Fix{\chi}=\Fix{\chi|_{\langle x,y \rangle}} *\langle
y_1,...,y_{n-2}\rangle$$ where each $y_i$ is equal to $a_i$ if $i
< r$ or $a_i g_i a_i^{-1}$ for some $g_i \in \langle
x,y,a_1,...,a_{i-1} \rangle$ otherwise. It is clear from the construction that the free factors $A_i$ are $\chi$-invariant, since we know that  
$K_i$ and $H_i$ are $\chi$-invariant. 

Now if $\chi|_{\langle x,y \rangle} $ is of polynomial growth,
then repeated applications of Lemma~\ref{poly} would show that
$\chi$ has polynomial growth. Therefore,  $\chi|_{\langle x,y
\rangle} $ is of exponential growth and has fixed subgroup of rank
1. Since automorphisms of the free group of rank 2 are geometric,
 $\chi|_{\langle x,y \rangle} $ is realised as a pseudo-Anasov
 automorphism of a once punctured torus. Thus we see that, up to a choice of basis,
 the fixed subgroup of  $\chi|_{\langle x,y \rangle} $ is $\langle
 [x,y] \rangle$. \endproof

The following shows how information about the fixed subgroup, when it is large, puts restrictions on the automorphism. 

\begin{lem}
\label{goodrtt} Let $\chi$ be an automorphism with fixed subgroup
as in Proposition~\ref{basis}. Then there exists a relative train
track representative of $\chi$ whose only exponential stratum is
the bottom stratum corresponding to (the conjugacy class of)
$\langle   x,y\rangle $. Moreover, this representative is
aperiodic (see \cite{Tits1} Definitions~3.1.7).
\end{lem}
\demo 
First as the smallest free factor
containing $xyx^{-1}y^{-1}$ is $\langle   x,y\rangle $, this must
be $\chi$-invariant, as $\chi$ sends free factors to free factors.
Now let $A_i=\langle   x,y,a_1...,a_i\rangle $ for $0 \leq i \leq
n-2$.

We claim that $A_i$ is $\chi$-invariant for all $i$. We have just shown
that this is true for $A_0$ and we argue by induction. Suppose the
claim holds for $A_i$. If $y_{i+1}=a_{i+1}$, then there is nothing
to prove. So assume that $a_{i+1} g_{i+1} a_{i+1}^{-1}$ which is
fixed by $\chi$. Note that $g_{i+1} \in A_i$ and by induction so
is $\chi(g_{i+1}$. Thus, $a_{i+1}^{-1} \chi(a_{i+1})$ is in the
normaliser of $A_i$, which is equal to $A_i$ as free factors are
malnormal. Hence $\chi(a_{i+1})= a_{i+1} u_{i+1}$ for some
$u_{i+1} \in A_i$. This proves the claim.

Now, by \cite{Tits1} Lemma~2.6.7, there is a relative train track representative of $\chi$, $f$ on $G$, such that 
$G$ has $f$-invariant subgraphs $G_1 \subseteq G_2 \subseteq \cdots G_{n-1}=G$ so that each $G_i$ is a connected subgraph whose fundamental group is equal to $A_i$, via the marking. 

Note that the difference in ranks between $A_{i+1}$ and $A_i$ is exactly one, so that each $G_{i+1}-G_i$ is a level stratum, for $i >1$. Note also that as $G_1$ has a rank $2$ fundamental group, $f_{G_1}$ must be of exponential growth, since otherwise $\phi$ would have polynomial growth. In particular, $f_{G_1}$ must be irreducible and also aperiodic. Hence the subgraphs $G_i$ give a complete stratification, where only the bottom stratum is exponential and corresponds to the conjugacy class of $\langle x,y \rangle$. Since this is aperiodic, the whole relative train track is aperiodic. \endproof

\begin{prop}
\label{invariant} Let $\chi, \phi \in \Aut{F_n}$ such that $\chi$
satisfies the hypotheses of Proposition~\ref{basis}. If $\Fix{\phi
\chi \phi^{-1}}=\Fix{\chi}$ then $\phi(\langle x,y\rangle )$ is a
conjugate of $\langle x,y\rangle $.
\end{prop}
\demo 
By Lemma~\ref{goodrtt}, $\chi$ has an aperiodic relative
train track representative with a unique exponential stratum at
the bottom. Thus $\chi$ has a unique attracting lamination, $\Lambda$, 
(see \cite{Tits1}, section 3), and since the bottom stratum of our relative train track corresponds
to $\langle x,y\rangle $, the attracting lamination of $\phi$ is
carried by the conjugacy class of $\langle x,y\rangle $. (We shall abuse notation slightly and call $\Lambda$ an attracting lamination for $\chi$ when it is strictly speaking an attracting lamination for the outer automorphism determined by $\chi$.)

Since we have only used Lemma~\ref{goodrtt} and the fixed
subgroup, we also deduce that the conjugacy class of $\langle
x,y\rangle $ also carries the unique attracting lamination, $\mathcal{L}$, of 
$\phi \chi \phi^{-1}$.

However, by inspection, the conjugacy class of $\phi(\langle
x,y\rangle )$ must also carry the lamination $\Phi(\Lambda)$, where $\Phi$ is the outer automorphism determined by $\phi$. However, $\Phi(\Lambda)$ is an attracting lamination for $\phi \chi \phi^{-1}$, since $\Phi$ induces a homeomorphism on $\mathcal{B}$ and sends bi-recurrent lines to bi-recurrent lines (see \cite{Tits1}, Definition~2.2.2 and Lemma~3.1.4). 

In other words $\Phi(\Lambda)=\mathcal{L}$ and $\phi(\langle x,y\rangle
)$ is conjugate to $\langle x,y\rangle $. \endproof

\section{Auto-fixed subgroups of rank $n-1$}

We shall need to quote one more result before we prove that the intersection of 1-auto-fixed subgroups is 1-auto-fixed if the rank of the intersection is at least $n-1$. 

\begin{thm}[\cite{MaVe}, Lemma~3.1]
\label{freefac} Let $\phi, \psi \in Aut(F_n)$. Then there exists a
positive integer $k$ such that $\Fix{\phi} \cap \Fix{\psi}$ is a
free factor of $\Fix{\phi \psi^k}$.
\end{thm}

\begin{thm}
\label{lessone} Let $\phi, \psi \in \Aut{F_n}$ and suppose that
the rank of $\Fix{\phi} \cap \Fix{\psi}$ is $n-1$. Then there
exists a $\chi \in \Aut{F_n}$ such that $\Fix{\chi}=\Fix{\phi}
\cap \Fix{\psi}$.
\end{thm}
\demo By Theorem~\ref{freefac} there exists an integer $k$ such that
$H=\Fix{\phi} \cap \Fix{\psi}$ is a free factor of $\Fix{\phi
\psi^k}$. Either $H=\Fix{\phi \psi^k}$ or $rank(\Fix{\phi
\psi^k})=n$. In the former case, we are done so we shall assume
that $\Fix{\phi \psi^k}$ has rank $n$ and so by Proposition~\ref{maxpoly} is UPG. 
Applying Theorem~\ref{freefac} again, there exists an
integer $m$ such that $H=\Fix{\phi \psi^k} \cap \Fix{\psi}$ is a
free factor of $\Fix{\psi (\phi \psi^k)^m}$. Thus, without loss of generality,
we may assume that $rank(\Fix{\phi})=rank(\Fix{\psi})=n$. Note that by \cite{DS}, as $\phi$ restricts to a finite order automorphism on $\Fix{\phi^k}$, $\Fix{\phi}$ is a free factor of $\Fix{\phi^k}$. But in our case, $rank(\Fix{\phi})=n$ and by Theorem~\ref{inertia}, $rank(\Fix{\phi^k}) \leq n$. Hence $\Fix{\phi^k}=\Fix{\phi}$ for all $k \neq 0$. By the same argument, $\Fix{\psi^k}=\Fix{\psi}$ for all $k \neq 0$.

If $\langle   \phi, \psi\rangle $ is a UPG subgroup then we are
done by Theorem~\ref{upg}. If every automorphism in $\langle   \phi, \psi\rangle $ has polynomial growth then, by 
\cite{Tits2} Proposition~3.5, $\langle   \phi, \psi\rangle $ has a finite index UPG subgroup. However, we would be done in that case since the fixed subgroups of
$\phi, \psi$ are unchanged by taking proper powers. Thus we may assume that there exists a $\chi \in \langle   \phi,
\psi\rangle $ of exponential growth. If $\Fix{\chi}=H$ then we are
done.


Now we can write $\Fix{\chi}=K*L$, where $K$ is an algebraic extension of $H$ (see \cite{KM} or \cite{V1}). 
But, as $H$ is inert, $rank(H) \leq rank(K) \leq rank(\Fix{\chi})$ and since $\chi$ has exponential growth, we know by Proposition~\ref{maxpoly}, that $rank(\Fix{\chi}) < n$. Thus $\Fix{\chi}=K$ is an algebraic extension of $H$. However, there are only finitely many of these, again by \cite{V1}, and since $H$ is fixed by both $\phi$ and $\psi$, these two automorphisms must permute the algebraic extensions of $H$. 
Thus $\Fix{\chi}$ is stabilised by some powers of $\phi$ and
$\psi$, which without loss we may assume to be $1$. (Note, in fact by
\cite{Tits2}, $\phi, \psi$ must always stabilise $\Fix{\chi}$.)

Hence, $\Fix{\phi \chi \phi^{-1}}=\Fix{\chi}=\Fix{\psi \chi \psi^{-1}}$.

By Proposition~\ref{basis}, there is a rank $2$ free factor $L$ of $F_n$ which is invariant under $\chi$ and such that the restriction of $\chi$ to $L$ is of exponential growth. By Proposition~\ref{invariant}, $\phi$ sends  $L $ to a conjugate of itself. By the same argument, 
$\psi$ also sends $L$ to a
conjugate. Let $\phi_0, \psi_0$ be automorphisms in the same outer
automorphism class as $\phi, \psi$ respectively which each leave
$L$ invariant.

By Theorem~\ref{res}, $\phi_0, \psi_0$ each restrict to a UPG
automorphism of $L$. Thus, up to a choice of
basis, $\phi_0$ induces an action of the
abelianisation of $L$ corresponding to the
matrix $\left(\begin{array}{cc} 1 & n
\\ 0 & 1
\end{array}\right)
$ for some integer $n$. With respect to this basis the induced
action of $\psi_0$ will correspond to some matrix
$\left(\begin{array}{cc} a & b \\ c & d
\end{array}\right),$ for some integers, $a,b,c,d$. As $\psi_0$ is also UPG,
the trace of this last matrix must equal $2$. 

By Theorem~\ref{freefac}, there exists a 
positive integer $k$ such that $H$ is a free factor of
$\Fix{\psi \phi^{k}}$. Either this automorphism has fixed subgroup exactly equal to $H$ and we are done, or it has fixed subgroup of rank $n$. In
the latter case, it will be UPG, by Proposition~\ref{maxpoly}. 
We shall show that this leads to a contradiction.

Note, again by \cite{Tits2} Proposition~4.6, $\psi_0 \phi_0^{k}$ induces a UPG automorphism on $L$. Hence the matrix $\left(\begin{array}{cc} a & b
\\ c & d
\end{array}\right). \left(\begin{array}{cc} 1 & n \\ 0 & 1
\end{array}\right)^{k}$ has trace $2$.

Evaluating this matrix we see that $\left(\begin{array}{cc} a & akn+ b \\
c & ckn+d
\end{array}\right)$ can only have trace $2$ if either $n=0$ or $c=0$. Note that the latter case forces $a=d=1$. 
Hence the image of $\langle   \phi_0,
\psi_0\rangle $ in $GL(2, \mathbb{Z})$, under the natural map,  is a subgroup consisting of
upper uni-triangular matrices. Thus, after choosing a basis $x,y$ for $L$, $\chi|_{L}$ is the automorphism  which sends $x $ to $g^{-1} x g$ and $y$ to $g^{-1} y x^t g$ for some integer, $t$, and some $g \in \langle x,y \rangle$. This contradicts the fact that $\chi|_{L} $ has exponential growth. (Here we are using the well known fact that the kernel of the map from $Aut(F_2)$ to $GL(2, \mathbb{Z})$ consists of precisely the inner automorphisms.) This contradiction completes the proof. 
\endproof

\begin{cor}
\label{rank3} Any auto-fixed subgroup of $F_3$ is $1$-auto-fixed.
\end{cor}
\demo By Thereom~\ref{finite}, it is enough to show that the intersection of two $1$-auto-fixed subgroups of $F_3$ is also $1$-auto-fixed. 

Let $\phi, \psi \in Aut(F_3)$ and let $H=\Fix{\phi} \cap \Fix{\psi}$. We need to show that $H=\Fix{\chi}$ for some $\chi \in Aut(F_3)$. By Theorem~\ref{inertia}, $rank(H)
\leq 3$. If $rank(H)=3$ then the result follows by \cite{MaVe}, Corollary~4.1. 
If $rank(H)=2$ then the result follows by \ref{lessone}. If
$rank(H)=1$ then $H=\langle   g\rangle $ and the inner automorphism
that is conjugation by $g$ will have fixed subgroup equal to $H$.
Finally, if $H$ is trivial then one can find automorphisms which
have no fixed points. An automorphism which cyclically permutes a
basis would be an example. 
\endproof

\Addresses
\recd

\end{document}

%% file: agtout.tex
\def\ifplaintex{\expandafter\ifx\csname documentclass\endcsname\relax}

\def\gtp{{\mathsurround=0pt\it $\cal G\mskip-2mu$eometry \&\ 
$\cal T\!\!$opology $\cal P\!$ublications}}  

\def\recd{{\small Received:\qua\receiveddate\ifx\reviseddate\relax
\else\qquad Revised:\qua\reviseddate\fi\par}} 

\def\lognumber#1{\def\thelognumber{#1}}
\def\volumenumber#1{\def\thevolumenumber{#1}}
\def\volumeyear#1{\def\thevolumeyear{#1}}
\def\papernumber#1{\def\thepapernumber{#1}}
\def\pagenumbers#1#2{\def\startpage{#1}\def\finishpage{#2}}
\def\published#1{\def\publishdate{#1}}

\def\received#1{\def\receiveddate{#1}}
\def\revised#1{\def\reviseddate{#1}}
\def\accepted#1{\def\accepteddate{#1}}

\let\\\par\let\thelognumber\relax\let\thevolumenumber\relax
\let\thepapernumber\relax\let\thevolumeyear\relax\let\startpage\relax
\let\finishpage\relax\let\publishdate\relax\let\receiveddate\relax
\let\reviseddate\relax\let\accepteddate\relax\let\theasciititle\relax
\let\theasciiauthors\relax
\let\theasciiabstract\relax

\let\theasciiemail\relax

\ifplaintex
\font\logobig=cmssbx10 scaled 3836
\font\logomed=cmssbx10 scaled 2557
\else
\font\logobig=cmssbx10 scaled 4200
\font\logomed=cmssbx10 scaled 2800
\fi

\long\def\makeagttitle{   
\count0=\startpage
\agt\hfill      
\hbox to 45truept{\vbox to 0pt{\vglue -13truept{\logomed A\kern -.37em{\logobig 
T}\kern -.38em G}\vss}\hss}
\break
{\small Volume \thevolumenumber\ (\thevolumeyear)
\startpage--\finishpage\nl
Published: \publishdate}

\vglue .25truein

{\parskip=0pt\leftskip 0pt plus
1fil\def\\{\par\smallskip}{\Large\bf\thetitle}\par\medskip} \vglue
0.05truein

{\parskip=0pt\leftskip 0pt plus 1fil\def\\{\par}{\sc\theauthors}
\par\medskip}%
 
\vglue 0.03truein 

{\small\leftskip 25truept\rightskip 25truept{\bf Abstract}\stdspace\theabstract

{\bf AMS Classification}\stdspace\theprimaryclass
\ifx\thesecondaryclass\relax\else; \thesecondaryclass\fi\par
{\bf Keywords}\stdspace \thekeywords\par}\vglue 7truept

}   

\ifplaintex
\font\phead=cmsl9 scaled 950
\font\pnum=cmbx10 scaled 913
\font\pfoot=cmsl9 scaled 950
\headline{\vbox to 0pt{\vskip -4.5mm\line{\small\phead\ifnum
\count0=\startpage ISSN 1472-2739 (on-line) 1472-2747 (printed)
\hfill {\pnum\folio}\else\ifodd\count0\def\\{ }%
\ifx\theshorttitle\relax\thetitle\else\theshorttitle\fi\hfill{\pnum\folio}
\else\def\\{ and }{\pnum\folio}\hfill\ifx\theshortauthors\relax\theauthors
\else\theshortauthors\fi\fi\fi}\vss}}
\footline{\vbox to 0pt{\vglue 0mm\line{\small\pfoot\ifnum\count0=\startpage
\copyright\ \gtp\hfill\else
\agt, Volume \thevolumenumber\ (\thevolumeyear)\hfill\fi}\vss}}
\else
\font=cmsl9 scaled 1050
\font\lnum=cmbx10 
\font=cmsl9 scaled 1050
\makeatletter
\def\@oddhead{{\small\lhead\ifnum\count0=\startpage ISSN 1472-2739 
(on-line) 1472-2747 (printed)\hfill {\lnum\number\count0}\else\ifodd\count0
\def\\{ }\ifx\theshorttitle\relax \thetitle \else\theshorttitle\fi\hfill
{\lnum\number\count0}\else\def\\{ and }{\lnum\number\count0}
\hfill\ifx\theshortauthors\relax 
\theauthors\else\theshortauthors\fi\fi\fi}}\def\@evenhead{\@oddhead}
\def\@oddfoot{\small\lfoot\ifnum\count0=\startpage\copyright\ \gtp\hfill\else
\agt, Volume \thevolumenumber\ (\thevolumeyear)\hfill\fi}
\def\@evenfoot{\@oddfoot}
\makeatother
\fi
\let\maketitlepage\makeagttitle
\let\makeshorttitle\maketitlepage
\let\maketitle\maketitlepage

\newwrite\gtoutfile
\long\gdef\makeheadfile{  
{\def\\{, }\def\s{ }
\immediate\openout\gtoutfile head.xxx
\immediate\write\gtoutfile{Proxy-for: \ifx\theasciiauthors\relax
\theauthors\else\theasciiauthors\fi\s<\ifx\theasciiemail\relax\theemail\else\theasciiemail\fi>}
\immediate\write\gtoutfile{\noexpand\\}
\immediate\write\gtoutfile{Authors: \ifx\theasciiauthors\relax
\theauthors\else\theasciiauthors\fi}
{\def\\{ }\immediate\write\gtoutfile{Title: \ifx\theasciititle\relax
\thetitle\else\theasciititle\fi}}
\immediate\write\gtoutfile{Subj-class: GT or SG, GR etc}
\immediate\write\gtoutfile{MSC-class: \theprimaryclass\ifx\thesecondaryclass\relax\else, \thesecondaryclass\fi}
\immediate\write\gtoutfile{Journal-ref: Algebraic and Geometric Topology \thevolumenumber\s
(\thevolumeyear) \startpage-\finishpage}
\immediate\write\gtoutfile{Comments: Published by Algebraic and
Geometric Topology at}
\immediate\write\gtoutfile{\s\s\s  http://www.maths.warwick.ac.uk/agt/AGTVol\thevolumenumber/agt-\thevolumenumber-\thepapernumber.abs.html}
\immediate\write\gtoutfile{\noexpand\\}
\immediate\write\gtoutfile{}
\ifx\theasciiabstract\relax
\immediate\write\gtoutfile{\theabstract}\else
\immediate\write\gtoutfile{\theasciiabstract}\fi
\immediate\write\gtoutfile{}
\immediate\write\gtoutfile{\noexpand\\}
\immediate\write\gtoutfile{}
\immediate\closeout\gtoutfile}}  

\def\maketitlepage{\makeagttitle\makeheadfile}
\let\makeshorttitle\maketitlepage
\let\maketitle\maketitlepage